\documentclass[12pt]{article}

\usepackage{amsmath}
\usepackage{amssymb}
\usepackage[mathscr]{eucal}
\usepackage{graphicx}
\usepackage{color}

\newtheorem{Theorem}{\sc Theorem}
\newtheorem{Definition}[Theorem]{\sc Definition}
\newtheorem{Proposition}[Theorem]{\sc Proposition}

\newtheorem{Corollary}[Theorem]{\sc Corollary}

\newtheorem{Example}[Theorem]{\sc Example}

\newcommand{\cQ}{\mbox{{${\cal Q}$}}}
\newcommand{\cM}{\mbox{{${\cal M}$}}}
\newcommand{\cP}{\mbox{{${\cal P}$}}}

\newcommand{\cS}{\mbox{{${\cal S}$}}}

\newcommand{\ve}{\mbox{{$\varepsilon$}}}
\newcommand{\R}{{\if mm {\rm I}\mkern -3mu{\rm R}\else \leavevmode
\hbox{I}\kern -.17em\hbox{R} \fi}}

\newcommand{\cT}{\mbox{{${\cal T}$}}}
\newcommand{\cC}{\mbox{{${\cal C}$}}}

\newcommand{\bx}{\mbox{\boldmath{$x$}}}

\def\sqr#1#2{{
    \vcenter{
         \vbox{\hrule height.#2pt
               \hbox{\vrule width.#2pt height#1pt \kern#1pt
                     \vrule width.#2pt
               }
               \hrule height.#2pt
         }
    }
}}

\def\bar{\overline}
\def\real{\mathbb{R}}

\def\lista#1
{{ \itemindent 0.0cm \labelsep .2cm \leftmargin 0.8cm \rightmargin
0.0cm \labelwidth 0.6cm \topsep 0.0mm
\parsep 0.0mm
\itemsep 0.0mm
\begin{list}{}
{ \setlength{\leftmargin}{.8cm} \setlength{\rightmargin}{0.0cm}
\setlength{\parsep}{0.0mm} \setlength{\topsep}{.0mm}
\setlength{\parskip}{.0cm} \setlength{\itemsep}{.0cm} }
{#1}\end{list}} }

\newcounter{theorem}

\begin{document}


\title{\bf On the Tykhonov Well-posedness of an Antiplane Shear Problem}

\vspace{22mm}
{\author{Mircea Sofonea$^{1}$\footnote{Corresponding author, E-mail : sofonea@univ-perp.fr}\, and\, Domingo A. Tarzia$^{2}$\\[6mm]
{\it \small $^1$ Laboratoire de Math\'ematiques et Physique}\\
{\it \small
	University of Perpignan Via Domitia}
\\{\it\small 52 Avenue Paul Alduy, 66860 Perpignan, France}		\\[6mm]		
{\it\small $^2$  Departamento de Matematica-CONICET}\\ {\it \small Universidad Austral}\\
		{\it \small Paraguay 1950, S2000FZF Rosario, Argentina}}

\date{}
\maketitle
\thispagestyle{empty}

\vskip 6mm

\noindent {\small{\bf Abstract.} We consider a boundary value problem which describes the  frictional antiplane shear of an elastic body. The process is static and friction is modeled with a slip-dependent version of Coulomb's law of dry friction. The weak formulation of the problem is in the form of a quasivariational inequality for the displacement field, denoted by  $\cP$. We associated to problem $\cP$ a boundary  optimal control problem, denoted by $\cQ$.  For  Problem $\cP$  we introduce the concept of well-posedness and for Problem
$\cQ$ we introduce the
concept of weakly and weakly generalized well-posedness, both associated to appropriate
Tykhonov triples. Our main result are Theorems \ref{t1} and \ref{t2}. Theorem \ref{t1} provides the well-posedness of Problem $\cP$ and, as a consequence, the continuous dependence of the solution with respect to the data. Theorem \ref{t2} provides the  weakly generalized well-posedness of Problem $\cQ$ and, under additional hypothesis, its weakly well posedness. The proofs of these theorems are based on arguments of compactness, lower semicontinuity, monotonicity and various estimates.
Moreover, we provide the mechanical interpretation of our well-posedness results.}

\vskip 6mm


\vskip 6mm


\vskip 15mm

\section{Introduction}\label{s1}
\setcounter{equation}0


The concept of Tykhonov well-posedness
was introduced in \cite{Ty} for a minimization problem and then
it has been generalized for different optimization problems, see for instance \cite{CKR, DZ, H, HY, L, Z},
The  well-posedness in the sense of  Tykhonov (well-posedness, for short) has been extended in the recent years to various mathematical problems like inequalities, inclusions, fixed point and saddle point problems. Thus, the well-posedness  of variational inequalities
was studied for the first time in \cite{LP1,LP2} and the study of well-posedness  of hemivariational inequalities was initiated in \cite{GM}.  References in the field include \cite{SX8,XHW}, among others.
A general framework which unifies the view  on the
well-posedness  for abstract problems in metric spaces  was recently introduced in \cite{SX9}. There,  the well-posedness concept has been introduced by using    families of  approximating sets $\{\Omega(\theta)\}_\theta$ indexed upon a positive parameter $\theta>0$, together with approximating sequences, which are sequences $\{u_n\}_n$ such that $u_n\in\Omega(\theta_n)$ for all $n\in\mathbb{N}$.

 Antiplane shear deformations are one of the simplest classes of
 deformations that solids can undergo: in antiplane shear (or
 longitudinal shear) of a cylindrical body, the displacement is
 parallel to the generators of the cylinder and is independent of
 the axial coordinate. We rarely actually load solids so as to cause them to deform
 in antiplane shear. However, the governing
 equations and boundary conditions for antiplane shear problems are
 beautifully simple and the solution will have many of the features
 of the more general case and may help us to solve the more complex
 problem too. For this reason, in the last years a considerable
 attention has been paid to the analysis of antiplane shear
 deformation, see for instance \cite{Ho1,Ho3}
 and the review article \cite{Ho2}. There,
 modern
 developments for the antiplane shear model and its applications
 are described for both linear and nonlinear solid materials,
 various physical settings (including dynamic effects) where the
 antiplane shear model shows promise for further development are
 outlined and some open questions concerning a challenging
 antiplane shear inverse problem in linear isotropic elastostatics
 are described. A reference in the study of antiplane contact problems with elastic and viscoelastic materials is the book \cite{SM-07}.

In this paper we consider a  mathematical model
which describes the antiplane shear of the cylinder in frictional contact. The model was already considered in \cite{SM-07} and, therefore we skip the mechanical assumptions used to derive the equations and boundary conditions involved. In a weak formulation, the model leads to an elliptic quasivariational inequality in which the unknown is the axial component of the displacement field. We associate to this problem a boundary optimal control problem.  Optimal control problems for variational inequalities  have been discussed in several works, including \cite{Ba,BT,F,Mig,MP, NST} and, more recently \cite{BT1,BT2,C}. Results on optimal control for antiplane contact problems with elastic materials can be found in \cite{MM1,MM3} and the references therein. There, the existence of optimal pairs was obtained, under various assumptions on the material's behaviour.

Our aim in this paper is two fold. The first one is to study the well-posedness of the quasivariational inequality above mentionned. The second one is to  study the well-posedness of the associated optimal control problem. In this way we complete the studies initiated in \cite{SM-07}, on one hand, and in \cite{MM1,MM3}, on the other hand.
Our main results are gather in Theorems \ref{t1} and \ref{t2}. The novelty of these results consists in the fact that we use the concept of Tykhonov triple in order to define the family of approximating sets and, in this way we provide the well-posedness of two relevant problems by using a functional framework which extends the  functional framework used in \cite{SX9}. As a consequence, we obtain new results, which represent a continuation  of the resuls obtained in \cite{MM1,MM3} and
give rise to interesting consequences and mechanical interpretation.

The rest of the manuscript is structured as follows. In Section \ref{s2} we introduce the antiplane shear problem, list the assumption on the data and state its variational formulation, denoted by $\cP$. Then, we associate to Problem $\cP$ a boundary optimal control problem, denoted by $\cQ$.
In Section \ref{s3}  we introduce
the concept of well-posedness for Problem $\cP$.
Moreover, we state and prove our first result, Theorem \ref{t1}, together with some consequences.
In Section \ref{s4}  we introduce
the concepts of weakly and weakly generalized well-posedness for Problem $\cQ$.
Then,  we state and prove our second result, Theorem \ref{t2}. The proofs of the theorems are based on arguments of monotonicity, compacteness, lower semicontinuity and various estimates. Finally, we
provide various mechanical interpretations of our well-posedness results.

\section{The antiplane shear model}\label{s2}
\setcounter{equation}0

We consider  the slip-dependent frictional antiplane shear model introduced in \cite[Section 9.2]{SM-07}. This model can be formulated as follows.

\medskip\noindent{\bf Problem} $\cM$. {\it
Find a displacement field $u:D\to \R$ such that}
	\begin{align}
	\label{I1} div\,(\mu\,\nabla u)+ f_0&=0\quad&\mbox{\rm in}\
	D, \\[2mm]
	\label{I2} u&=0\quad &\mbox{\rm on}\ \Gamma_1,\\[2mm]
	\label{I3} \mu\,\partial_\nu u&=f_2\quad&\mbox{\rm on}\
	\Gamma_2,\\[2mm]
	\label{I4} |\mu\,\partial_\nu\,u|\leq
	g(|u|),\quad\mu\,\partial_\nu u&=-g(|u|)\,\frac{u}{|u|}\quad
	&\mbox{\rm if}\ u\neq 0\quad \mbox{\rm on}\ \Gamma_3.
	\end{align}

Recall that Problem $\cM$ describes the equilibrium of an elastic cylinder of transversal section $D$
under the action of axial  body forces of density $f_0$ and surface traction of density $f_2$. Here
$D$ is assumed to be a regular domain in $\R^2$ with smooth boundary
$\partial D=\Gamma$, composed of
three measurable sets $\overline{\Gamma}_1$, $\overline{\Gamma}_2$ and
$\overline{\Gamma}_3$, with the mutually disjoint relatively open
sets $\Gamma_1$, $\Gamma_2$ and $\Gamma_3$, such that the
one-dimensional measure of $\Gamma_1$, denoted by $meas\,\Gamma_1$,
is strictly positive.
Equation (\ref{I1}) represents the equation of equilibrium in which $\mu$ denotes the Lam\'e coefficient, (\ref{I2}) represents the displacement boundary condition and (\ref{I3}) represents the traction
boundary condition. Finally, (\ref{I4}) is a static version of slip dependent Coulomb's
law in which $g$ is a positive function, the friction bound, assumed to
depend on the slip, i.e., $g=g(|u|)$. Note that here, in (\ref{I1})--(\ref{I4}),  and  below in this paper,  we skip the dependence of various functions with respect to the spatial variable $\bx\in D\cup\Gamma$.

For the variational analysis of Problem $\cM$ we use standard notation for Lebesque and Sobolev spaces. We use the symbols $``\to"$ and $"\rightharpoonup"$ to indicate the strong and weak convergence in various spaces which will be indicated below. Moreover, we use the space $V$ given by
\begin{equation}\label{SpV}
V=\{\,v\in H^1(D)\,:\, v=0\,\mbox{ on }\Gamma_1\,\}
\end{equation}
and we denote by $(\cdot,\cdot)_V$, the inner product induced on $V$ by the inner product of $H^1(D)$. In addtion, $\|\cdot\|_V$  and $0_V$ will represent the associated norm and the zero element of $V$. It is well known that $(V,\|\cdot\|_V)$ is a real Hilbert space. Moreover, by the
Friedrichs-Poincar\'e inequality and the Sobolev trave theorem we have
\begin{eqnarray}
&&\label{FP}
	\|v\|_{V}\leq c_0\,\|\nabla
	v\|_{L^2(D)^2}, \\ [2mm]
&&\label{trace}
\|v\|_{L^2(\Gamma_3)}\le c_3\,\|v\|_V
\end{eqnarray}
for all $v\in V$, respectively. Here
$c_0$ and $c_3$ are positive constants which depend on $D$, $\Gamma_1$ and $\Gamma_3$.

We now list the assumption of the data. First,  we assume that the Lam\'e coefficient and the densities of the body forces and surface traction satisfy the following conditions.
\begin{eqnarray}
&&\label{mu}\mu\in L^\infty(D)\mbox{ and there exists }\mu^*>0 \mbox{
	such that }\label{I10}\\
&&\qquad\mu(\bx)\geq \mu^* \mbox{ a.e. }\bx\in D.\nonumber\\ [2mm]
&&\label{f0}
f_0\in L^2(D).\\ [2mm]
&&\label{f2} f_2\in L^2(\Gamma_2).
\end{eqnarray}
Moreover, we  assume that the  friction bound
is such that
\begin{equation}\label{g}
\left\{
\begin{array}{l}
(a)\ g:\Gamma_3\times\mathbb{R}\to\mathbb{R}_+.\\[2mm]
(b)\ \mbox{There exists }L_g>0\mbox{ such that }\\[0mm]
\qquad|g(\bx,r_1)-g(\bx,r_2)|\leq L_g\,|r_1-r_2|\\[0mm]
\qquad\quad\forall\, r_1,\,r_2\in \R,\,\mbox{ a.e. }\bx\in \Gamma_3.\\[2mm]
(c)\ \mbox{The mapping }\bx\mapsto g(\bx,\,r)\\
\qquad\mbox{ is Lebesgue measurable on }\Gamma_3,\,\forall\,r\in
\mathbb{R}.\\[2mm] (d)\ \mbox{The mapping }\bx\mapsto g(\bx,0)\mbox{
	belongs to }L^2(\Gamma_3).
\end{array}
\right.
\end{equation}

\noindent
Finally, we assume that the following smallness condition hold:
\begin{equation}\label{sm}
L_gc_0^2c_3^2<\mu^*.
\end{equation}

Define the bilinear form  $a:V\times V\to \R$, the function $j:V\times V\to\R$ and the element $f\in V$ by equalities
\begin{eqnarray}
&&\label{a} a(u,v)=\int_{D}\,\mu\,\nabla u\cdot\nabla
v\,dx\qquad\forall\,u,\,v\in V,\\[2mm]
&&\label{j} j(u,\,v)=\int_{\Gamma_3}g(|u|)\,|v|\,da\qquad\forall\, u,\,v\in V,\\ [2mm]
&&\label{f} (f,v)_V=\int_{D}f_0v\,dx+\int_{\Gamma_2}f_2v\,da \qquad\forall\, v\in V,
\end{eqnarray}
Then, using standard arguments, it is easy to
derive the following variational formulation of Problem $\cM$.

\medskip\noindent{\bf Problem} $\cP$. {\it
	Find a displacement field $u\in V$ such that}
	\begin{equation}\label{iv}
	a(u,v-u)+j(u,v)-j(u,u)\geq (f,v-u)_V\quad\forall\, v\in V.
	\end{equation}
	
\medskip	
A function $u\in V$ which satisfies (\ref{iv})
and is called a {\it weak solution} to Problem ${\cal M}$.
Note that the function $j$  in (\ref{iv}) depends on the solution $u$ and, for this reason, we refer to  (\ref{iv})  as a quasivariational inequality.

Now, in order to formulate a boundary optimal control for Problem $\cP$, we assume that $\mu$, $f_0$ and $g$ are given and, for any $f_2$, we use notation (\ref{f}). Then, the set of admissible pairs for inequality (\ref{iv}) is  given by
\begin{equation}\label{2}
{\cal V}_{ad} = \{\,(u, f_2)\in V\times L^2(\Gamma_2) \ \mbox{such that}\   (\ref{iv})\  \mbox{holds}\,\}.
\end{equation}
We also consider the cost function ${\cal L}$ given by
\begin{eqnarray}
&&\label{L}\hspace{-12mm}{\cal L}(u,f_2)=a_0\int_{D}|u-\phi|^2\,da+a_2\int_{\Gamma_2}|f_2|^2\,da=a_0\|u-\phi\|^2_{L^2(D)}+a_2\|f_2\|^2_{L^2(\Gamma_2)}
\end{eqnarray}
for all  $u\in V$, $f_2\in L^2(\Gamma_2)$. Here $\phi$ is a given element in $V$ and  $a_0$, $a_2$ are strictly positive constants, i.e.,
\begin{eqnarray}
&&\label{y1}\varphi\in V,\\ [2mm]
&&\label{y2}a_0,\ a_2>0.
\end{eqnarray}

With the notation above we consider the following optimal control problem.

\medskip\medskip\noindent
{\bf Problem} ${\cal Q}$.  {\it Find $(u^*, f^*_2)\in {\cal V}_{ad}$ such that}
\begin{equation}\label{3}
{\cal L}(u^*,f^*_2)=\min_{(u,f_2)\in {\cal V}_{ad}} {\cal L}(u,f_2).
\end{equation}

With this choice for ${\cal L}$, the mechanical interpretation of  Problem $\cQ$ is the following:
given a frictional antiplane shear  described by the quasivariational inequality  (\ref{iv}) with the  data $\mu$, $f_0$ and $g$ which satisfy condition (\ref{mu}), (\ref{f0}), (\ref{g}), respectively,
we are looking for a surface traction $f_2\in L^2(\Gamma_2)$ such that
the displacement $u$ of the cylinder solution  is as close as possible, to the ``desired  normal displacement" $\phi$. Furthermore, this choice has to fulfill a minimum expenditure condition which is taken into account by the last  term in (\ref{L}). In fact, a compromise policy between the two aims (``$u$ close to $\phi$" and ``minimal data $f_2$") has to be found and the relative importance of each criterion with respect to the other is expressed by the choice of the weight coefficients $a_0,\, a_2> 0$.

\medskip

The unique solvability of Problem $\cP$ as well as the solvability of Problem $\cQ$ represent a consequence of assumptions (\ref{mu})--(\ref{sm}), (\ref{y1}), (\ref{y2}) and can be obtained by using standard arguments, already used in \cite{SM-07} and \cite{MM1,MM3}, respectively.
Our aim in what follows is to study the well-posedness of these problems in the sense of Tykhonov and to derive some consequences together with their mechanical interpretation. To this end we need to prescribe three ingredients: a set of indices, a familly of approximating sets and a convergence for the sequences of indices. These three ingredients lead us to introduce the new  concept of Tykhonov triple, both for Problem $\cP$ and Problem $\cQ$, in Sections \ref{s3} and \ref{s4}, respectively.

\section{A well-posedness result}\label{s3}
\setcounter{equation}0

Everywhere in this section we assume that $(\ref{mu})$--$(\ref{sm})$ hold. Moreover, we note that the concepts of
Tykhonov triple, approximating sequence and well-posedness we introduce below in this section refer to Problem $\cP$, even if we do not mention it explicitly.
Extending by our previous work  \cite{SX9}, we consider the following definitions.

\begin{Definition}\label{d0}
A Tykhonov triple for Problem $\cP$ is a mathematical object of the form $\cT=(I,\Omega,\cC) $ where  $I$ is a given set, $\Omega:I\to 2^V$ and $\cC\subset\cS(I)$.	
\end{Definition}

Note that in this definition  and below in this paper $\cS(I)$ represents the set of sequences of elements of $I$ and  $2^V$ denotes the set of parts of the space $V$. A typical emelent of $I$ will be denoted by $\theta$ and a typical element of $\cS(I)$ will be denoted by $\{\theta_n\}_n$. For any $\theta\in I$ we refer to the set $\Omega(\theta)\subset V$ as  an approximating set and $\cC$ representes the so-called convergence criterion.

\begin{Definition}\label{d1}	
Given a  Tykhonov triple $\cT=(I,\Omega,\cC) $, 	a sequence $\{u_n\}_n\subset V$ is called an approximating sequence  if there exists a sequence  $\{\theta_n\}_n\subset\cC$ such that   	 $u_n\in \Omega(\theta_n)$ for each $n\in\mathbb{N}$.
\end{Definition}

\begin{Definition}\label{d2} The Problem $\cP$  is said to be (strongly) well-posedness with respect to  the Tykhonov triple $\cT=(I,\Omega,\cC) $ if it has a unique solution  and every approximating sequence  converges  in $V$ to this solution.
\end{Definition}

Note that the concept of  approximating sequence above depends on the Tykhonov triple $\cT$ considered and, for this reason, we sometimes  refer to approximating sequences with respect to $\cT$. As a result, the concept of  well-posedness depends on the Tykhonov triple $\cT$.  We also remark that the choice of this triple is crucial for the analysis of the well-posedness of Problem $\cP$.  In what follows we construct an exemple of such triple, used in the study of Problem $\cP$.

\begin{Example}\label{ex1}
Take $\cT=(I,\Omega,\cC)$ where
\begin{eqnarray}
\label{I}
&&\hspace{-10mm}I=\{\,\theta=(\ve,\widetilde{f_0},\widetilde{f_2},\widetilde{g})\ :\ \ve\ge 0,\ \widetilde{f_0},\, \widetilde{f_2},\ \widetilde{g}\ \mbox{\rm satisfy\ (\ref{f0}),  (\ref{f2}) and (\ref{g})}\,\},\\ [2mm]
&&\label{O}\hspace{-10mm}\Omega(\theta)=\{\, u\in V\ :\   \label{ive}
a(u,v-u)+\widetilde{j}(u,v)-\widetilde{j}(u,u)+\ve\|u\|_V\|v-u\|_V\\ [0mm]
&&\hspace{28mm}\geq (\widetilde{f},v-u)_V\quad \forall\, v\in V\,\}\nonumber
\end{eqnarray}
where, for a given $\theta=(\ve,\widetilde{f_0},\widetilde{f_2},\widetilde{g})\in I$, $\widetilde{j}$ and $\widetilde{f}$ are defined  by $(\ref{j})$, $(\ref{f})$, replacing  $g$ with $\widetilde{g}$, $f_0$ with $\widetilde{f_0}$ and $f_2$ with $\widetilde{f}_{2}$.
Moreover, by definition, a sequence $\{\theta_n\}_n\subset \cS(I)$ with   $\theta_n=(\ve_n,f_{0n},f_{2n},g_n)$
belongs to $\cC$ if and only if the following hold:
\begin{eqnarray}
&&\label{C1}
\ve_n\to 0,\\ [2mm]
&&\label{C2} f_{0n}\rightharpoonup f_0\ \ {\rm in}\ L^2(D),\\ [2mm]
&&\label{C3}
 f_{2n}\rightharpoonup f_2\ \ {\rm in}\ L^2(\Gamma_2)
\end{eqnarray}
as  $n\to\infty$ and, in addition,
\begin{eqnarray}
&&\label{C4}
\left\{
\begin{array}{ll}
\mbox{\rm for any}\ n\in\mathbb{N}\ \mbox{there exists }\alpha_n\ge 0\ \mbox {\rm and}\ \beta_n\ge 0\ \mbox{ such
	that
}\\[2mm]
\mbox{\rm (a) }
|g_{n}(x,r)-g(x,r)|
\leq
\alpha_n+ \beta_n|r|
\quad\forall\, r\in\R,\ {\rm a.e.}\  \bx\in \Gamma_3.\\[2mm]
\mbox{\rm (b) }\displaystyle\lim_{n\to\infty}\alpha_n=\lim_{n\to\infty}\beta_n=0.
\end{array}
\right.
\end{eqnarray}
\end{Example}

\medskip
Our main result in this section is the following.

\begin{Theorem}\label{t1} Assume that $(\ref{mu})$--$(\ref{sm})$ hold. Then Problem $\cP$  is  well-posed
with respect to the Tykhonov triple  $\cT$  in Example $\ref{ex1}$.
\end{Theorem}	

\noindent	
{\it Proof.} Following Definition \ref{d2},  the proof is carried out in two main steps.

\medskip\noindent i) {\it Unique solvability of Problem $\cP$.} First, we use assumption (\ref{mu}) and inequality (\ref{FP}) to see that the bilinear form $a$ is continuous and coercive, i.e.,
\begin{equation}\label{ba}
a(v,v)\ge \frac{\mu^*}{c_0^2}\|v\|_V^2\qquad\forall\, v\in V.
\end{equation}
Next, we use assumption (\ref{g}) and inequality (\ref{trace}) to see that for all $\eta\in V$	$ j(\eta,\cdot):V\to\R$ is a continuous seminorm and, moreover,
\begin{eqnarray}
&&\label{bj} j(\eta_1,v_2)-j(\eta_1,v_1)+j(\eta_2,v_1)-j(\eta_2,v_2)\\ [2mm] &&\qquad \le L_gc_3^2\,\|\eta_1-\eta_2\|_V\|v_1-v_2\|_V\quad
\forall \,\eta_1,\,\eta_2,\,v_1,\,v_2\in V.\nonumber
\end{eqnarray}
In addition, recall the smallness assumption (\ref{sm}).
The existence of a unique solution $u\in V$ to Problem $\cP$ follows now by using a standard argument of quasivariational inequalities, see, for instance \cite{SM-07}.

\medskip\noindent ii) {\it Convergence of approximating sequences.}
To proceed, we
consider  an approximating sequence for the Problem $\cP$, denoted  $\{u_n\}_n$. Then, according to Definition \ref{d1}  it follows that
there exists a sequence $\{\theta_n\}_n$ of elements of $I$, with  the general term denoted by $\theta_n=(\ve_n,f_{0n},f_{2n},g_n)$, such that $u_n\in\Omega(\theta_n)$
for each $n\in\mathbb{N}$ and, moreover, (\ref{C1})--(\ref{C4})
hold.
Our aim in what follows is to
prove the convergence
\begin{equation}\label{m}
u_n\to u\qquad{\rm in}\quad V,\quad{\rm as}\quad n\to\infty.
\end{equation}
To this end we proceed in three intermediate steps that we present below.

\medskip\noindent {\rm ii-a)} {\it A boundedness result.} Let $n\in\mathbb{N}$ be fixed.
Then, exploiting definition (\ref{O}) we deduce that
\begin{eqnarray}\label{ivn}
&&a(u_n,v-u_n)+j_n(u_n,v)-j_n(u_n,u_n)+\ve_n\|u_n\|_V\|v-u_n\|_V\\[2mm]
&&\qquad\qquad\geq (f_n,v-u_n)_V\qquad \forall\, v\in V\nonumber
\end{eqnarray}
where
\begin{eqnarray}
&&\label{jn} j_n(u,\,v)=\int_{\Gamma_3}g_n(|u|)\,|v|\,da\qquad\forall\, u,\,v\in V,\\ [2mm]
&&\label{fn} (f_n,v)_V=\int_{D}f_{0n}v\,dx+\int_{\Gamma_2}f_{2n}v\,da \qquad\forall\, v\in V.
\end{eqnarray}

 We take $v=0_V$ in (\ref{ivn}) to obtain that
\[a(u_n,u_n)+j_n(u_n,u_n)\le\ve_n\|u_n\|_V^2+(f_n,u_n)_V.\]
Using now inequalities (\ref{ba}),  $j_n(u_n,u_n)\ge 0$ and $(f_n,u_n)_V\le \|f_n\|_V \|u_n\|_V$
yields
\[\Big(\frac{\mu^*}{c_0^2}-\ve_n\Big)\|u_n\|_V\le \|f_n\|_V.\]
Finally, we use the convergences (\ref{C1})--(\ref{C3}) to deduce that the sequence $\{u_n\}_n$ is bounded in $V$, i.e., there exists $k$ which does not depend on $n$ such that
\begin{equation}\label{bu}
\|u_n\|_V\le k.
\end{equation}

\medskip\noindent {\rm ii-b)} {\it Weak convergence.}
The bound (\ref{bu}) allows us to deduce that there  exists an element $\widetilde{u}\in V$ and a subsequence of $\{u_n\}_n$, again denoted $\{u_n\}_n$, such that
\begin{equation}\label{9}
{u}_n\rightharpoonup\widetilde{u}\quad \mbox{ \rm{in} }\
V\quad \mbox{ \rm{as} }\ n\to\infty.
\end{equation}
Let $n\in\mathbb{N}$ and $v\in V$.  Then,  (\ref{ivn}) implies that
\begin{eqnarray}\label{ip}
&&a(u_n,v)+j_n(u_n,v)-j_n(u_n,u_n)+\ve_n\|u_n\|_V\|v-u_n\|_V\\[2mm]
&&\qquad\qquad\geq (f_n,v-u_n)_V+a(u_n,u_n).\nonumber
\end{eqnarray}
We now use the assumption (\ref{C4}), the convergence (\ref{9}) and the compactness of the trace operator to see that
\begin{equation*}
g_n(|u_n|)\to g(|\widetilde{u}|),\quad|u_n|\to |\widetilde{u}|\quad  {\rm in}\quad L^2(\Gamma_3) \quad{\rm as}\quad n\to\infty
\end{equation*}
and, therefore,
\begin{equation}\label{z1}
j_n(u_n,v)-j_n(u_n,u_n)\to j(\widetilde{u},v)-j(\widetilde{u},\widetilde{u})\quad{\rm as}\quad n\to\infty.
\end{equation}
Moreover, the convergences (\ref{C1}) and (\ref{9}) imply that
\begin{equation}\label{z2}
\ve_n\|u_n\|_V\|v-u_n\|_V\to 0\quad{\rm as}\quad n\to\infty.
\end{equation} Finally, the  convergences (\ref{C2}), (\ref{C3}) and (\ref{9}) combined with standard compactness arguments imply that
\begin{equation}\label{z3}
(f_n,v-u_n)_V\to (f,v-\widetilde{u})_V\quad{\rm as}\quad n\to\infty.
\end{equation}

On the other hand,  using the properties of the form $a$ we have that
\begin{equation}\label{z3n}
a(u_n,v)\to a(\widetilde{u},v)\quad{\rm as}\quad n\to\infty,\quad\forall\, v\in V.
\end{equation}
Moreover, since $a(u_n-\widetilde{u},u_n-\widetilde{u})\ge 0$ we deduce that \[a(u_n,u_n)\ge a(\widetilde{u},u_n)+a(u_n,\widetilde{u})-a(\widetilde{u},\widetilde{u}).\]
Using now (\ref{z3n})  and the convergence (\ref{9}) we find that
\begin{equation}\label{z4}
\liminf_{n\to\infty}\, a(u_n,u_n)\ge a(\widetilde{u},\widetilde{u}).
\end{equation}
We now pass to the lower limit in inequality (\ref{ip}) and use (\ref{z1})--(\ref{z4}) to deduce that
\begin{equation}\label{11}a(\widetilde u,v-\widetilde u)+j(\widetilde u,v)-j(\widetilde u,\widetilde u) \ge(f,v-\widetilde u)_V.
\end{equation}
Next, it follows from (\ref{11}) that $\widetilde{u}$ is a solution of inequality (\ref{iv})
and, by the uniqueness of the solution of this inequality, we obtain that
\begin{equation}
\label{eut} \widetilde{u}=u.
\end{equation}

A careful analysis, based on the arguments above, reveals that $u$ is the weak limit of any weakly convergent subsequence of the sequence  $\{u_n\}_n$.
Therefore, using a standard argument we deduce that
the whole sequence
$\{u_n\}_n$ converges weakly in $V$ to $u$ as
$n\to \infty$, i.e.,
\begin{equation}\label{9n}
{u}_n\rightharpoonup{u}\quad \mbox{ \rm{in} }\
V\quad \mbox{ \rm{as} }\ n\to\infty.
\end{equation}

\medskip\noindent
{\rm ii-c)} {\it Strong convergence.}	
Let $n\in\mathbb{N}$ be given. We take $v=u$ in inequality
(\ref{ivn}) to see that
\begin{equation}\label{3v} a(u_n,u_n-u)\le j_n(u_n,u)-j_n(u_n,
u_n)+\ve_n\|u_n\|_V\|u-u_n\|_V+(f_n,u_n-u)_V.
\end{equation}
Next, we use inequatities (\ref{ba})  and  (\ref{3v}) to find that
\begin{eqnarray*}
	&&\label{yw}\frac{\mu^*}{c_0^2}\,\|u_n-u\|_V^2\leq a(u_n-u,u_n-u)=a(u_n,u_n-u)-a(u,
	u_n-u)_V\nonumber\\[2mm]
	&&\quad\leq j_n(u_n,u)-j_n(u_n,
	u_n)+\ve_n\|u_n\|_V\|u-u_n\|_V+(f_n,u_n-u)_V-a(u,
	u_n-u).\nonumber
\end{eqnarray*}
We now pass to the upper limit in this inequality and use
the convergences (\ref{z1})--(\ref{z4}), (\ref{9n}) and equality  (\ref{eut}) to deduce that
\begin{equation*}
\|u_n-u\|_V\to 0\quad\mbox{
	\rm{as} }\ n\to \infty.
\end{equation*}
This convergence concludes
the proof.\hfill$\Box$

\medskip
Theorem \ref{t1} provides an useful tool in the study of the continuous dependence of the solution with respect to the data. Various exemples can be considered but, for simplicity, we restrict ourselves to present two examples which make the object of the corrolaries below.

\begin{Corollary}	\label{cor1} Assume that $(\ref{mu})$--$(\ref{sm})$ hold and, for each $n\in\mathbb{N}$, assume that $f_{0n}\in L^2(D)$, $f_{2n}\in L^2(\Gamma_3)$ and $g_n$ are given, $g_n$ satisfies condition  $(\ref{g})$ and, moreover, $(\ref{C2})$, $(\ref{C3})$ and $(\ref{C4})$ hold. Then
the solution of Problem $\cP$ with the data  $f_{0n}$, $f_{2n}$ and $g_n$, denoted by $u_n$, converge to the solution of  Problem $\cP$, i.e., $(\ref{m})$ holds.
\end{Corollary}

\noindent	
{\it Proof.} We use the Thkhonov triple $\cT=(I,\Omega,\cC)$ in Example \ref{ex1}. For each $n\in\mathbb{N}$ denote $\theta_n=(0,f_{0n},f_{2n},g_n)\in I$. Then it is easy to see that the sequence $\{\theta_n\}_n$ belongs to the set $\cC$. Moreover,  by  definition, it follows that for each $n\in\mathbb{N}$ the element $u_n\in V$ satisfies the quasivariational inequality
\begin{eqnarray*}
a(u_n,v-u_n)+j_n(u_n,v)-j_n(u_n,u_n)\geq (f_n,v-u_n)_V\quad\forall\, v\in V
 \end{eqnarray*}
 where, recall, $j_n$ is defined by (\ref{jn}) and $f_n$ is defined by (\ref{fn}) . We now use  (\ref{O}) to deduce that $u_n\in\Omega(\theta_n)$ for each $n\in\mathbb{N}$ and, hence, Definition \ref{d1} shows that $\{u_n\}_n$ is an approximating sequence for Problem $\cP$. Using now Theorem \ref{t1} and Definition \ref{d2} we deduce that the convergence  $(\ref{m})$ holds, which concludes the proof.\hfill$\Box$

\begin{Corollary}	\label{cor2} Assume that $(\ref{mu})$--$(\ref{sm})$ hold and, for each $n\in\mathbb{N}$, assume that $\mu_n$ is given, satisfies condition $(\ref{mu})$ and, moreover,
\begin{equation}\label{cmu}
\mu_n\to\mu\quad \mbox{ \rm{in} }\
L^\infty(D)\quad \mbox{ \rm{as} }\ n\to\infty.
\end{equation}Then
the solution of Problem $\cP$ with the data  $\mu_n$, denoted $u_n$, converge to the solution of  Problem $\cP$, i.e., $(\ref{m})$ holds.
\end{Corollary}

\noindent	
{\it Proof.} Let $n\in\mathbb{N}$ and $v\in V$ be given. Then
\begin{eqnarray}\label{q1}
a_n(u_n,v-u_n)+j(u_n,v)-j(u_n,u_n)\geq (f,v-u_n)_V\quad\forall\, v\in V
\end{eqnarray}
where $a_n$ is defined by equality (\ref{a}) in which $\mu$ is replaced by $\mu_n$, that is
\begin{equation}\label{q2}a_n(u,v)=\int_{D}\,\mu_n\,\nabla u\cdot\nabla
v\,dx\qquad\forall\,u,\,v\in V.
\end{equation}
We now combine (\ref{q1}) and (\ref{q2}) to see that
\begin{eqnarray*}
&&a(u_n,v-u_n)+\int_{D}\,(\mu_n-\mu)\,\nabla u_n\cdot\nabla
(v-u_n)\,dx+j(u_n,v)-j(u_n,u_n)\geq (f,v-u_n)_V
\end{eqnarray*}
and, since
\begin{eqnarray*}
&&\int_{D}\,(\mu_n-\mu)\,\nabla u_n\cdot\nabla
(v-u_n)\,dx\\ [2mm]
&&\quad\le \|\mu_n-\mu\|_{L^\infty(\Gamma_3)}\|\nabla u_n\|_{L^2(D)^2} \|\nabla (v-u_n)\|_{L^2(D)^2}\\ [3mm]
&&\qquad\le  \|\mu_n-\mu\|_{L^\infty(D)}\|u_n\|_V\|v-u_n\|_V,
\end{eqnarray*}
we deduce that
\begin{eqnarray}\label{q3}
	&&a(u_n,v-u_n)+j(u_n,v)-j(u_n,u_n)+\ve_n\|u_n\|_V\|v-u_n\|_V\geq (f,v-u_n)_V
\end{eqnarray}
where
\begin{equation}\label{q4}
\ve_n=\|\mu_n-\mu\|_{L^\infty(D)}.
\end{equation}

 We use the Tykhonov triple $\cT=(I,\Omega,\cC)$ in Example \ref{ex1}. For each $n\in\mathbb{N}$ denote $\theta_n=(\ve_n,f_0,f_{2},g)\in I$. Then, using (\ref{q4}) and (\ref{cmu})  it follows that
 (\ref{C1}) holds and, moreover, the sequence $\{\theta_n\}_n$ belongs to the set $\cC$. In addtion, inequality (\ref{q3})
  and  (\ref{O}) show that $u_n\in\Omega(\theta_n)$ for each $n\in\mathbb{N}$. Therefore, using Definition \ref{d1} we deduce that $\{u_n\}_n$ is an approximating sequence for Problem $\cP$. Finally,$  $ using Theorem \ref{t1} and Definition \ref{d2} we deduce that the convergence  $(\ref{m})$ holds, which concludes the proof. \hfill$\Box$

\medskip
In additional to the mathematical interest in the convergence  result in Corollary \ref{cor1} it is important from mechanical point of view since show that small perturbations on the density of body forces and surface tractions together with small perturbation on the friction bound give rise to small perturbations of the weak solution of the antiplane contact problem $\cM$. A similar comment can be made on Corollary  \ref{cor2} which shows the
continuous dependence of the weak solution of Problem $\cM$ with respect to the Lam\'e coefficient $\mu$.
We conclude from here the importance of the well-posedness concept in the sense of Tykhonov in the study of antiplane contact problems.

\medskip
We proceed with an elementary example in which an analytic expression of the solution of Problem $\cP$  can be given.

\begin{Example}\label{ex1n}
Consider the one-dimensional version of Problem $\cM$ in which $D=(0,1)$, $\Gamma_1=\{0\}, \Gamma_2=\emptyset$, $\Gamma_3=\{1\}$, $\mu>0$, $f_0\in \real$ and $g>0$. With these data the problem consists to find a function $u:[0,1]$ such that
\begin{eqnarray}\label{sof}
-\mu u''=f_0,\quad\  u(0)=0,\quad\  |u'(1)|\le g,\quad\  u'(1)=-g\,\frac{u(1)}{|u(1)|}\ \ {\rm if}\ \ u(1)\ne 0.
\end{eqnarray}
where, here and below, the prime represents the derivative with respect to the statil variable $x$. A simple calculus show that the solution of the nonlinear boundary value problem $(\ref{sof})$ is given by	
\begin{equation}
u(x)=\left\{\begin{array}{ll}
-\frac{f_0}{2\mu}\,x^2+\Big(\frac{f_0}{\mu}+g\Big)x\qquad{\rm if}\quad f_0<-2\mu g\\[3mm]-\frac{f_0}{2\mu}\,x^2+\frac{f_0}{2\mu}\,x\hskip19.5mm{\rm if}\quad -2\mu g\le f_0\le 2\mu g,\\[3mm]
-\frac{f_0}{2\mu}\,x^2+\Big(\frac{f_0}{\mu}-g\Big)x\qquad{\rm if}\quad f_0>2\mu g
\end{array}\right.\label{uu}
\end{equation}
Moreover, note that this solution is the solution of the corresponding Problem $\cP$ which, in this particular case, is stated as follows: find $u\in V$ such that
\begin{eqnarray*}
	\mu\int_0^1 u'(v'-u')\,dx+g(|v(1)|)-g(|u(1)|)\ge\int_0^1f_0(v-u)\,dx\qquad\forall\,v\in V
\end{eqnarray*}
where $V=\{\,v\in H^1(0,1)\,:\, v(0)=0\,\}$.
It is easy to see that the function $(\ref{uu})$ depends continuously on the data $\mu$, $f_0$ and $g$, which represents a validation of the convergence results in Corollaries $\ref{cor1}$ and $\ref{cor2}$.

\end{Example}

We end this section with an example which illustrate the sensitivity of the well-posedness concept of Problem $\cP$ with respect to the choice of the Tykhonov
triple.

\begin{Example}\label{ex7}
	Let $\bar{f}_0\in L^2(D)$ be given and  take $\cT=(I,\Omega,\cC)$ where
	\begin{eqnarray}
	\label{I7}
	&&\hspace{-10mm}I=L^2(\Omega),\\ [2mm]
	&&\label{O7}\hspace{-10mm}\Omega(\theta)=\{\, u\in V\ :\   \label{ivep}
	a(u,v-u)+j(u,v)-j(u,u)\geq (\widetilde{f},v-u)_V\quad \forall\, v\in V\,\}
	\end{eqnarray}
	where, for a given $\theta=\widetilde{f_0}\in I$,  $\widetilde{f}$ is defined  by $(\ref{f})$, replacing  $f_0$ with $\widetilde{f_0}$.
	Moreover,
	\begin{equation}
	\hspace{-18mm}\label{C27}
	\cC=\{\, \{f_{0n}\}_n\subset\cS(L^2(\Omega))\ :\  f_{0n}\rightharpoonup \bar{f}_0\ \ {\rm in}\ L^2(D)\quad{\rm as }\quad n\to\infty.\,\}
	\end{equation}
\end{Example}

\medskip
We have the following. result

\begin{Proposition}\label{t1n} Assume that $(\ref{mu})$--$(\ref{sm})$ hold.
Then Problem $\cP$  is  well-posed
	with respect to the Tykhonov triple  $\cT$  in Example $\ref{ex7}$ if and only if $f_0=\bar{f}_0$.
\end{Proposition}	

\noindent	
{\it Proof.} Assume that $f_0=\bar{f}_0$
and let $\{u_n\}_n$ be an approximating sequence for Problem $\cP$ with respect to the Tykhonov triple $\cT$ in Example \ref{ex7}. Then, there exists a sequence $\{f_{0n}\}_n\subset L^2(D)$ such that $f_{0n}\rightharpoonup {f}_0 \ \ {\rm in}\ L^2(D)\ {\rm as }\ n\to\infty$ and, moreover, for each $n\in\mathbb{N}$, $u_n$ satisfies inequality
\begin{eqnarray}\label{q17}
a(u_n,v-u_n)+j(u_n,v)-j(u_n,u_n)\geq (f_n,v-u_n)_V\quad\forall\, v\in V,
\end{eqnarray}
in which $f_n$ is defined by (\ref{f}), replacing $f_0$ with $f_{0n}$.
It follows from here that $\{u_n\}_n$ is an approximating sequence for Problem $\cP$ with respect to the Tykhonov triple $\cT$ in Example \ref{ex1}, corresponding to the sequence $\{\theta_n\}_n$ with $\theta_n=(0,f_{0n},f_2,g)$ for each $n\in\mathbb{N}$. Therefore, using Theorem \ref{t1} and Definition \ref{d2} we deduce that $u_n\to u$ in $V$ as $n\to \infty$, where $u$ denotes the solution of $\cP$. This shows that Problem $\cP$ is well posed with the Tykhonov triple (\ref{I7})--(\ref{C27}).

Conversely, assume that Problem $\cP$ is well posed with the Tykhonov triple in Example \ref{ex7}
and denote by $\bar{\cP}$ the quasivariational inequality  in Problem $\cP$ obtained by replacing $f_0$ with $\bar{f}_0$.

Let $\{u_n\}_n$ be an approximating sequence for Problem $\cP$ with respect to the Tykhonov triple in Example \ref{ex7}.
Then, since $\cP$ is well posed with $\cT$, we deduce that
\begin{equation}\label{za1}
u_n\to u\qquad{\rm in}\qquad V\quad{\rm as}\quad n\to\infty
\end{equation} where, recall, $u$ denotes the solution to Problem $\cP$.
On the other hand, by Definition \ref{d1},
we know that there exists a sequence $\{f_{0n}\}_n\subset L^2(D)$ such that $f_{0n}\rightharpoonup \bar{f}_0 \ \ {\rm in}\ L^2(D)\ {\rm as }\ n\to\infty$ and, moreover, for each $n\in\mathbb{N}$, $u_n$ satisfies inequality (\ref{q17}) where,
again, $f_n$ is defined by (\ref{f}), replacing $f_0$ with $f_{0n}$. Next, for each $n\in\mathbb{N}$ denote $\theta_n=(0,f_{0n},f_2,g)\in I$. Then, from (\ref{q17}) it follows that $\{u_n\}_n$ is an approximating sequence from the Problem $\bar{\cP}$ and, therefore, using Theorem \ref{t1} we deduce that
\begin{equation}\label{za2}
u_n\to \bar{u}\qquad{\rm in}\qquad V\quad{\rm as}\quad n\to\infty
\end{equation}
where  $\bar{u}$ denotes the solution to Problem $\bar{\cP}$.

We now combine (\ref{za1}) and (\ref{za2}) to deduce that $\bar{u}=u$ Moreover, we recall that
$u$ denotes the solution to Problem $\cP$ with the data $\mu$, $f_0$, $f_2$ and $g$ while
$\bar{u}$ is the solution to Problem $\cP$ with the data $\mu$, $\bar{f}_0$, $f_2$ and $g$. Therefore, besides (\ref{iv}), $u$ satisfies the inequality
\begin{equation}\label{za4}
a(u,v-u)+j(u,v)-j(u,u)\geq (\bar{f},v-u)_V\quad\forall\, v\in V
\end{equation}
where
\begin{equation}\label{za5} (\bar{f},v)_V=\int_{D}\bar{f}_0v\,dx+\int_{\Gamma_2}f_2v\,da \qquad\forall\, v\in V.
\end{equation}
Let $w\in C_0^\infty(\Omega)$ be given. We take $v=u+w$ in (\ref{iv}) then $v=u-w$ in (\ref{za4})
and add the resulting inequalities to obtain that
\begin{equation}\label{za6}
j(u,u+w)-j(u,u)+j(u,u-w)-j(u,u)\ge(f-\bar{f},w)_V
\end{equation} Next, since $w$ vanishes on $\Gamma_3$, it follows that $j(u,u+w)=j(u,u-w)=j(u,u)$ and, therefore, (\ref{za6}) yields $(f-\bar{f},w)_V\le 0.$ A similar argument shows that
$(f-\bar{f},w)_V\ge 0$. We conclude from here that
$(f-\bar{f},w)_V=0$ and, using definitions (\ref{f}), (\ref{za5}) it follows that
\[\int_D(f_0-\bar{f}_0)w\,dx=0.\]
We now use a density argument to find that
$f_0-\bar{f}_0$ which concludes the proof.
\hfill$\Box$

\section{Weakly well-posedness results }\label{s4}
\setcounter{equation}0

Everywhere in this section we assume that $(\ref{mu})$, $(\ref{f0})$, $(\ref{g})$,  $(\ref{sm})$, $(\ref{y1})$ and $(\ref{y2})$ hold. Moreover, we underlie that the concept of
Tykhonov triple, approximating sequence and generalized well-posedness we introduce below in this section refear to Problem $\cQ$, even if we do not mention it explicitly.

\begin{Definition}\label{d0n}
	A Tykhonov triple for Problem $\cQ$ is a mathematical object of the form $\cT=(I,\Omega,\cC) $ where  $I$ is a given set, $\Omega:I\to 2^{V\times L^2(\Gamma_2)}$ and $\cC\subset\cS(I)$.	
\end{Definition}

Recall that in this definition  $\cS(I)$ represents the set of sequences of elements of $I$ and  $2^{V\times L^2(\Gamma_2)}$ denotes the set of parts of $V\times L^2(\Gamma_2)$. For any $\theta\in I$ we refer to the set $\Omega(\theta)\subset V\times L^2(\Gamma_2)$ as  an approximating set and $\cC$ representes the so-called convergence criterion.

\begin{Definition}\label{d1n}	
	Given a  Tykhonov triple $\cT=(I,\Omega,\cC) $, 	a sequence $\{(u_n,f_{2n})\}_n\subset {V\times L^2(\Gamma_2)}$ is called an approximating sequence if there exists a sequence  $\{\theta_n\}_n\subset\cC$, such that   	$(u_n,f_{2n})\in \Omega(\theta_n)$ for each $n\in\mathbb{N}$.
	
\end{Definition}

\begin{Definition}\label{d3n} The Problem $\cQ$  is said to be weakly  well-posedness with respect to  the Tykhonov triple $\cT=(I,\Omega,\cC) $ if it has a unique solution  and every approximating sequence converges weakly in $V\times L^2(\Gamma_2)$ to this solution.
\end{Definition}

\begin{Definition}\label{d2n} The Problem $\cQ$  is said to be weakly generalized well-posedness with respect to  the Tykhonov triple $\cT=(I,\Omega,\cC) $ if it has at least one  solution  and every approximating sequence contains a subsequence which converges weakly in $V\times L^2(\Gamma_2)$ to some point of the solution  set.
\end{Definition}

As in the previous section, we remark that the concept of  approximating sequence above depends on the Tykhonov triple  $\cT$ and, for this reason, we sometimes we refer to approximating sequences with respect $\cT$. As a result, the concept of  weakly and weakly generalized well-posedness depend on the Tykhonov triple $\cT$.  In what follows we construct a relevant exemple of such triple that we use int the study of Problem $\cQ$.

\begin{Example}\label{ex2}
	
Take $\cT=(I,\Omega,\cC)$ where
\begin{eqnarray}
\label{Ic}
&&\hspace{-10mm}I=\{\,\theta=(\ve,\widetilde{f_0},\widetilde{g},\widetilde{\phi})\ :\ \ve\ge 0,\ \widetilde{f_0},\ \widetilde{g},\ \widetilde{\phi} \ \mbox{\rm satisfy\ (\ref{f0}), (\ref{g}) and (\ref{y1})}\,\},\\ [2mm]
&&\label{Oc}\hspace{-10mm}\Omega(\theta)=\{\, (u^*,f_2^*)\in  {\cal V}_{ad}(\theta)\ :\
{\cal L}_\theta(u^*,f^*_2)=\min_{(u,f_2)\in {\cal V}_{ad}(\theta)} {\cal L}_\theta(u,f_2)\,\}
\end{eqnarray}
where, for a given $\theta=(\ve,\widetilde{f_0},\widetilde{g}, \widetilde{\phi})\in I$, the cost functional ${\cal L}_\theta$ and the set ${\cal V}_{ad}(\theta)$ and are defined as follows:
\begin{eqnarray}
&&\label{w0}
{\cal L}_\theta(u,f_2)=a_0\|u-\widetilde{\phi}\|^2_{L^2(D)}+a_2\|f_2\|^2_{L^2(\Gamma_2)},\\ [3mm]
&&\label{w1}(u,f_2)\in{\cal V}_{ad}(\theta) \ \Longleftrightarrow \quad u\in V,\quad f_2\in L^2(\Gamma_2)\ \ {\rm and}\\ [1mm]
&&
a(u,v-u)+\widetilde{j}(u,v)-\widetilde{j}(u,u)+\ve\|u\|_V\|v-u\|_V\geq (\widetilde{f},v-u)_V\quad \forall\, v\in V.\nonumber
\end{eqnarray}
Here, $\widetilde{j}$ and $\widetilde{f}$ are defined  by $(\ref{j})$, $(\ref{f})$, replacing  $g$ with $\widetilde{g}$ and $f_0$ with $\widetilde{f_0}$.
Finally, by definition, a sequence $\{\theta_n\}_n\subset \cS(I)$ with   $\theta_n=(\ve_n,f_{0n},g_n,\phi_n)$
belongs to $\cC$ if and only if $(\ref{C1})$, $(\ref{C2})$ and $(\ref{C4})$ hold and, moreover,
\begin{equation}\label{C5}
\phi_n\rightharpoonup\phi\ \ {\rm in}\ V.
\end{equation}

\end{Example}

\medskip
Our main result in this section is the following.

\begin{Theorem}\label{t2} Assume that $(\ref{mu})$--$(\ref{sm})$, $(\ref{y1})$ and $(\ref{y2})$ hold. Then Problem $\cQ$  is weakly  generalized  well-posed
with respect to the Tykhonov triple in Example $\ref{ex2}$.	
\end{Theorem}

\noindent	
{\it Proof.}  Following Definition \ref{d2n} the proof is carried out in two main steps.

\medskip\noindent i) {\it Solvability of Problem $\cQ$.} The existence of solutions to
Problem $\cQ$ follows from standard arguments that we resume in the following. Let
\begin{equation}\label{omega}
\omega=\inf_{(u,f_2)\in {\cal V}_{ad}} {\cal L}(u, f)\in\mathbb{R}
\end{equation}
and let $\{(u_n,f_{2n})\}_n\subset {\cal V}_{ad}$ be a minimizing sequence for the functional  ${\cal L}$, i.e.,
\begin{equation}\label{n}
\lim_{n\to\infty}\,{\cal L}(u_n, f_{2n})=\omega.
\end{equation}

We claim that the sequence $\{f_{2n}\}_n$ is bounded in $L^2(\Gamma_2)$.
Arguing by contradiction, assume  that $\{f_{2n}\}_n$ is not bounded in $L^2(\Gamma_2)$. Then, passing to a subsequence still denoted $\{f_{2n}\}_n$, we have
\begin{equation}\label{8}
\|f_{2n}\|_{L^2(\Gamma_2)}\to +\infty \quad\text{as}\quad n\to +\infty
\end{equation}and, using  (\ref{L}) it turns that
\begin{equation}\label{mm}
\lim_{n\to +\infty}{\cal L}(u_n, f_n)= +\infty.
\end{equation}
Equalities (\ref{n}) and (\ref{mm}) imply that $\omega=+\infty$ which is in contradiction with (\ref{omega}).

We conclude from above that the sequence $\{f_{2n}\}_n$ is bounded in $L^2(\Gamma_2)$, as clamed. Therefore, we deduce that there exists $f^*_2\in L^2(\Gamma_2)$ such that, passing to a subsequence still denoted $\{f_{2n}\}$, we have
\begin{equation}\label{ww}
f_{2n}\rightharpoonup f^*_2\quad\text{in}\quad L^2(\Gamma_3)\quad\text{as}\quad n\to +\infty.
\end{equation}

Let $u^*$ be the solution of the quasivariational inequality (\ref{iv}) for $f_2=f^*_2$.
Then,  by the definition (\ref{2}) of the set ${\cal V}_{ad}$ we have
\begin{equation}\label{sol-opt}
(u^*, f^*_2)\in {\cal V}_{ad}.
\end{equation}
Moreover, using  the convergence (\ref{ww}) and Corollary \ref{cor1}  it follows that
\begin{equation}\label{www}
u_n \to u^*\quad\text{in}\quad V\quad \text{as}\quad n\to +\infty.
\end{equation}
We now use the convergences (\ref{ww}), (\ref{www}) and the weakly lower semicontinuity of the functional ${\cal L}$, guaranteed by the fact that ${\cal L}$ is a convex continuous functional, to deduce that
\begin{equation}\label{xxx}
\liminf\limits_{n\to +\infty}{\cal L}(u_n,f_{2n})\geq{\cal L}(u^*,f^*_2).
\end{equation}
It follows from (\ref{n}) and (\ref{xxx}) that
\begin{equation}\label{xy}
\omega\geq{\cal L}(u^*,f^*_2).
\end{equation}
In addition, (\ref{sol-opt})  and (\ref{omega}) yield
\begin{equation}\label{xz}
\omega\le{\cal L}(u^*,f^*_2).
\end{equation}
We now combine inequalities  (\ref{xy})  and (\ref{xz}) to see that (\ref{3}) holds, which
shows that $(u^*,f_2)$ is a solution to Problem $\cQ$.

\medskip\noindent ii) {\it Convergence of approximating sequences.}
To proceed, we
consider  an approximating sequence for the Problem $\cQ$, denoted  by $\{(u_n^*, f_{2n}^*)\}_n$. Then, according to Definition \ref{d1n}  it follows that
there exists a sequence $\{\theta_n\}_n$ of elements of $I$, with $\theta_n=(\ve_n,f_{0n},g_n,\phi_n)$,  such that $(u_n,f_{2n})\in\Omega(\theta_n)$
for each $n\in\mathbb{N}$ and, moreover, (\ref{C1}),
(\ref{C2}), (\ref{C4}) and (\ref{C5})
hold.  In order to simplify the notation, for each $n\in\mathbb{N}$ we write ${\cal V}^n_{ad}$ and
${\cal L}_n$ instead of  ${\cal V}^n_{ad}(\theta_n)$ and ${\cal L}_{\theta_n}$, respectively.  Then, exploiting the definition (\ref{Oc}) we deduce that  $(u^*_n,f_{2n}^*)\in {\cal V}_{ad}^n$ and
\begin{eqnarray}
&&\label{zc}
{\cal L}_n(u^*_n,f^*_{2n})\le {\cal L}_n(u_n,f_{2n})
\end{eqnarray}
for each couple of functions $(u_n,f_{2n})\in  {\cal V}_{ad}^n$, i.e., for each couple of functions $(u_n,f_{2n})\in V\times L^2(\Gamma_2)$ which satisfies inequality (\ref{ivn}) in which, recall, $j_n$ and $f_n$ are defined by (\ref{jn}) and (\ref{fn}), respectively, and for each $n\in\mathbb{N}$.

We shall prove that there exists a subsequence of the sequence $\{(u_{n}^{*}, f_{2n}^{*})\}_n$, again denoted by $\{(u_{n}^{*}, f_{2n}^{*})\}_n$, and an element   $(u^*,f^*_2)\in V\times L^2(\Gamma_3)$ such that
\begin{eqnarray}
&&\label{se1}
f_{2n}^* \rightharpoonup f^*_2\quad\mbox{\rm in}\quad L^2(\Gamma_3)\quad\mbox{\rm as}\quad n\to \infty, \\[2mm]
&&\label{se2}
u_n^*\rightarrow u^*\quad \mbox{\rm in}\quad V\quad\mbox{\rm as}\quad n\to \infty, \\[2mm]
&&\label{se3}
(u^*,f^*_2)\quad {\rm is\ a\ solution\ of\ Problem}\ {\cal Q}.
\end{eqnarray}
To this end we proceed in four intermediate steps that we present below.

\medskip\noindent {\rm ii-a)} {\it A boundedness result.} We claim that the sequence  $\{f_{2n}^*\}_n$ is bounded in $L^2(\Gamma_2)$.
Arguing by contradiction, assume  that $\{f_n^*\}_n$ is not bounded in $L^2(\Gamma_2)$. Then, passing to a subsequence still denoted $\{f_n^*\}_n$, we have
\begin{equation}\label{8n}
\|f_{2n}^*\|_{L^2(\Gamma_2)}\to +\infty\quad\text{as}\quad n\to +\infty.
\end{equation}
We  use definition (\ref{w1}) of the cost functional ${\cal L}_n={\cal L}_{\theta_n}$ to deduce that
\begin{equation}\label{20}
\lim_{n\to \infty}{\cal L}_n(u_n^*, f_{2n}^*)= +\infty.
\end{equation}
Next, let $\bar{u}_n$ be the solution of the variational inequality
\begin{eqnarray}\label{ipz}
\bar{u}_n\in V,\quad a(\bar{u}_n,v-\bar{u}_n)+j_n(\bar{u}_n,v)-j_n(\bar{u}_n,\bar{u}_n)\geq (\bar{f}_n,v-\bar{u}_n)_V\quad\forall\, v\in V
\end{eqnarray}
in which
\begin{eqnarray}
&&\label{fnz} (\bar{f}_n,v)=\int_{D}f_{0n}v\,dx+\int_{\Gamma_3}f_{2}v\,da \qquad\forall\, v\in V.
\end{eqnarray}
Note that Theorem \ref{t1} guarantees that this solution exists and is unique, for each $n\in\mathbb{N}$. Moreover,  using (\ref{C2}) and (\ref{C4}) it follows that the sequence $\{\bar{u}_n\}_n$ is an approximating sequence from Problem $\cP$, with respect to the Tykhonov triple in Exemple \ref{ex1}, corresponding to the sequence $\{\theta_n\}_n$ with $\theta_n=(0,f_{0n},f_2,g_n)$ for all $n\in\mathbb{N}$. Therefore, using Theorem \ref{t1} it follows that
\begin{equation}
\bar{u}_n\rightarrow u\quad \mbox{\rm in}\quad V\quad\mbox{\rm as}\quad n\to \infty
\end{equation}
and, using  (\ref{C5}) and the definition of   the functional ${\cal L}_n$ we deduce that
\begin{equation}\label{e5}
{\cal L}_n(\bar{u}_n, f_2)\to{\cal L}(u, f_2) \quad\mbox{\rm as}\quad n\to \infty.
\end{equation}

On the other hand (\ref{ipz})  and (\ref{w1}) imply that the pair $(\bar{u}_n,f_2)$ satisfies inequality  (\ref{ivn}) with $f_{2n}=f_2$, i.e., $(\bar{u}_n,f_2)\in {\cal V}^n_{ad}$. Therefore (\ref{zc}) implies that
\begin{eqnarray}
&&\label{e6}
{\cal L}_n(u^*_n,f^*_{2n})\le {\cal L}_n(\bar{u}_n,f_{2})\qquad\forall\, n\in\mathbb{N}.
\end{eqnarray}
We now pass to the limit in (\ref{e6}) as $n\to\infty$ and use  (\ref{20}) and
(\ref{e5}) to obtain a contradiction.

\medskip\noindent {\rm ii-b)} {\it Two convergence results.} We conclude from step ii-a) that
the sequence $\{f_{2n}^*\}_n$ is bounded in $L^2(\Gamma_2)$ and, therefore,  we can find	a subsequence again denoted $\{f_{2n}^*\}_n$ and an element $f^*_2\in L^2(\Gamma_2)$ such that (\ref{se1}) holds.
	Denote by $u^*$ the solution of Problem ${\cal P}$ for $f_2=f^*_2$ and note that definition (\ref{2}) implies that
	\begin{equation}\label{se6}
	(u^*, f^*_2)\in {\cal V}_{ad}.
	\end{equation}
	Moreover, assumptions (\ref{C2}), (\ref{C4}) and the convergence (\ref{se1})   show that $\{u^*_n\}_n$ is an approximating sequence
	with the Tykhonov triple in Exemple \ref{ex1}, corresponding to the sequence $\{\theta_n\}_n$ with $\theta_n=(0,f_{0n},f^*_{2n},g_n)$ for all $n\in\mathbb{N}$.  Therefore,
	the well-\-posedness of Problem $\cP$, guaranteed by Theorem \ref{t1}, imply that (\ref{se2}) holds, too.

	\medskip\noindent {\rm ii-c)} {\it Optimality of the limit.}
	We now prove that $(u^*,f^*_2)$ is a solution to the optimal control problem ${\cal Q}$.
	To this end we use the convergences (\ref{se1}), (\ref{se2}), (\ref{C5}) and the weakly lower semicontinuity of the  functional
	$z\to \|z\|^2_{L^2(\Gamma_3)}$ on the space $L^2(\Gamma_2)$ to see that
	\begin{equation}\label{se8}
	{\cal L}(u^*,f^*_2)\leq\liminf_{n\rightarrow\infty}{\cal L}_n(u_n^*,f_{2n}^*).
	\end{equation}
	Next, we fix a  solution $({u}^*_0,{f}^*_{02})$ of Problem ${\cal Q}$ and,
	in addition, for each $n\in\mathbb{N}$  we denote by ${u}_n^0$ be the unique element of $V$ which satisfies the inequality
	\begin{eqnarray}\label{z23}
	&&
	a(u_n^0,v-u_n)+j_n(u_n^0,v)-j_n(u_n^0,u_n^0)\geq ({f}^0_n,v-u_n)_V\quad \forall\, v\in V,\nonumber
	\end{eqnarray}
in which
\begin{eqnarray*}
 (f_n^0,v)_V=\int_{D}f_{0n}v\,dx+\int_{\Gamma_3}f^*_{02}v\,da \qquad\forall\, v\in V.
\end{eqnarray*}	
It follows from here that $u_n^0$  satisfies inequality
\begin{eqnarray}\label{z24}
&&
a(u_n^0,v-u_n)+j_n(u_n^0,v)-j_n(u_n^0,u_n^0)+\ve_n\|u_n^0\|_V\|v-u^0\|_V\nonumber\\ [0mm]
&&\hspace{12mm}\geq ({f}^0_n,v-u_n)_V\quad \forall\, v\in V,\nonumber
\end{eqnarray}
too. Therefore, $({u}_n^0,{f}^*_{02})\in {\cal V}_{ad}^n$ and, using the optimality of the pair
	$(u_n^*,f_{n2}^*)$, (\ref{zc}), we find that
	\begin{equation*}
	{\cal L}_n(u_n^*,f_{2n}^*)\leq{\cal L}_n({u}_n^0,{f}^*_{02})\qquad\forall\, n\in\mathbb{N}.
	\end{equation*}
	We pass to the upper limit in this inequality to see that
	\begin{equation}\label{se9}
	\limsup_{n\rightarrow\infty}{\cal L}_n(u_n^*, f_{2n}^*)\leq \limsup_{n\rightarrow\infty}{\cal L}_n({u}_n^0,{f}^*_{02}).
	\end{equation}

	Next, remark that $\{u^0_n\}_n$ is an approximating sequence for Problem $\cP$ with respect to the Tykhonov triple in Example \ref{ex1} corresponding to the sequence $\{\theta_n\}_n$ with $\theta_n=(0,f_{0n},f_{02}^*,g_n)$ for all $n\in\mathbb{N}$. Therefore, Theorem \ref{t1} guarantees that
	\begin{equation*}
	{u}_n^0 \to {u}^*_0\quad\text{in}\quad V\quad\text{as}\quad n\to\infty.
	\end{equation*}
	Using now this convergence and assumption (\ref{C5}) yields
	\begin{equation}\label{se10}
	\lim_{n\rightarrow\infty}{\cal L}_n({u}_n^0,{f}^*_{02})={\cal L}({u}^*_0,f^*_{02}).
	\end{equation}
	We now use  (\ref{se8})--(\ref{se10}) to see that
	\begin{equation}\label{se10n}
	{\cal L}(u^*,f^*_2)\leq {\cal L}({u}^*_0,{f}^*_{02}).
	\end{equation}
	
	On the other hand, since $({u}^*_0,{f}^*_{02})$ is a solution of Problem ${\cal Q}$,  we have
	\begin{equation}\label{3p}
	{\cal L}({u}^*_0,{f}^*_{02})=\min_{(u,f_2)\in {\cal V}_{ad}} {\cal L}(u,f_2).
	\end{equation} and, therefore,
	inclusion (\ref{se6})
	implies that
	\begin{equation}\label{se11}
	{\cal L}({u}^*_0,{f}^*_{02})\le {\cal L}(u^*,f^*_2).
	\end{equation}
	We now combine the inequalities  (\ref{se10n}) and (\ref{se11}) to see that
	\begin{equation}\label{se16}
	{\cal L}(u^*, f^*_2)={\cal L}({u}^*_0,{f}^*_{02}).
	\end{equation}
	Next, we use relations (\ref{se6}), (\ref{se16}) and (\ref{3p}) to see that (\ref{se3}) holds.
	
	\medskip\noindent {\rm ii-d)} {\it End of proof.} We remark that the convergences (\ref{se1})  and (\ref{se2})  imply the weak convergence (in the product Hilbert space $V\times L^2(\Gamma_2)$) of the sequence $\{(u_n^*,f_{2n}^*)\}_n$ to the element $(u^*,f_2^*)$. Theorem \ref{t2} is now a direct consequence of Definition \ref{d2n}. \hfill$\Box$

\medskip
A direct consequence of Theorem \ref{t2} is the following.

\begin{Corollary}	\label{cor3}
	
Assume that $(\ref{mu})$--$(\ref{sm})$, $(\ref{y1})$, $(\ref{y2})$  hold and, moreover, assume that Problem $\cQ$ has a unique solution. Then Problem $\cQ$  is weakly  well-posed
with respect to the Tykhonov triple in Example $\ref{ex2}$.		
\end{Corollary}

\noindent	
{\it Proof.} Let $(u^*,f_2^*)$ be the unique solution to Problem $\cQ$ and let $\{(u_n^*, f_{2n}^*)\}_n$ be
an approximating sequence for the Problem $\cQ$, with respect to the Tyhoonov triple in Example \ref{ex2}. First, it follows from the proof of Theorem \ref{t2} that the sequence  $\{f_{2n}^*\}_n$ is bounded in $L^2(\Gamma_2)$. Therefore, using arguments similar to those used in step i) of the proof of Theorem \ref{t1}   we deduce that the sequence $\{u_{n}^*\}_n$ is bounded in $V$. We conclude from here that the sequence $\{(u_n^*, f_{2n}^*)\}_n$ is bounded in the product space $V\times L^2(\Gamma_2)$.
Second, a careful analysis of the proof of Theorem \ref{t2} reveals that $(u^*,f_2^*)$ is the weak limit  (in $V\times L^2(\Gamma_2)$) of any weakly convergent subsequence of the sequence   $\{(u_n^*, f_{2n}^*)\}_n$. The two properties above allow us to  use a standard argument in order to deduce   that
the whole sequence
 $\{(u_n^*, f_{2n}^*)\}_n$ converges weakly in $V\times L^2(\Gamma_2)$ to  $(u^*,f_2^*)$, as
$n\to \infty$. Corollary \ref{cor3} is now a direct consequence of Definition \ref{d3n}.\hfill$\Box$

\medskip
We now proceed with the following convergence result.

\begin{Corollary}	\label{cor4} Assume that $(\ref{mu})$--$(\ref{sm})$, $(\ref{y1})$, $(\ref{y2})$  hold and,  moreover, assume that for each $\widetilde{\phi}\in V$, Problem $\cQ$ has a unique solution. In addition, for each $n\in\mathbb{N}$, assume that $\phi_n$ is given and the convergence  and $(\ref{C4})$ holds. Then
	the solution of Problem $\cQ$ with the data  $\phi_n$, denoted by $(u_n^*,f_{2n}^*)$, converges weakly to the solution of  Problem $\cQ$, in the space $V\times L^2(\Gamma_2)$.
\end{Corollary}

\noindent	
{\it Proof.}
We use the Tykhonov triple $\cT=(I,\Omega,\cC)$ in Example \ref{ex2}, again. For each $n\in\mathbb{N}$ denote $\theta_n=(0,f_0,g,\phi_n)\in I$. Then, using (\ref{C5}) it is easy to see that the sequence $\{\theta_n\}_n$ belongs to the set $\cC$. Moreover, note that in this case ${\cal V}^n_{ad}={\cal V}_{ad}$, for all $n\in\mathbb{N}$. In addition, using the statement of Problem $\cQ$ we see that
\[(u_n^*,f_{2n}^*)\in V_{ad}, \quad{\rm and}\quad
{\cal L}_n(u_n^*,f_{2n}^*)\le {\cal L}_n(u_n,f_{2n})
\quad\forall\,(u_n,f_{2n})\in V_{ad}.\]
It follows now from the definition (\ref{Oc}) that
$(u_n^*,f_{2n}^*)\in\Omega(\theta_n)$ for each $n\in\mathbb{N}$ and, therefore,  Definition \ref{d1n} shows that $\{(u_n^*,f_{2n}^*)\}_n$ is an approximating sequence for Problem $\cQ$. We  now use Corollary \ref{cor4} and Definition \ref{d3n} to conclude the proof.\hfill$\Box$

\medskip
We end this section with the following comments, remarks and mechanical interpretation of Theorem \ref{t2} and Corollaries \ref{cor3} and \ref{cor4}.

\medskip
{\rm (i)} First,  recall that the convergence result (\ref{se2}) represents the strong convergence of the sequence $\{u_n\}_n$ to the element $u$, in the space $V$. Therefore, Theorem \ref{t2} provides more than the weakly generalized well-posedness of Problem $\cQ$ with respect to the Tykhonov triple in Example \ref{ex2}. Indeed, according to  Definition \ref{d2n}, to obtain the
weakly generalized well-posedness of Problem $\cQ$
we need only the weak convergence  $u_n\rightharpoonup u$ in $V$ as $n\to\infty$, which is obviously
implied by the strong convergence (\ref{se2}). A similar comment can be made concerning Corollaries \ref{cor3} and \ref{cor4}.

\medskip
{\rm (ii)} Second, recall that in general Problem $\cQ$ does not have a unique solution. The reason arises in the fact that
 the optimal control $\cQ$ is equivalent to  the problem of finding $u^*\in V$ and $f_2^*\in L^2(\Gamma_3)$ such that
\begin{equation}
\label{Jn}
u^*=u(f_2^*)\quad {\rm and}\quad J(f_2^*)=\min_{f_2\in L^2(\Gamma_2)} J(f_2),
\end{equation}
where  $u(f_2)$ represents the solution of Problem $\cP$ with the data $\mu$, $f_0$, $f_2$, $g$ and
$J:L^2(\Gamma_2)\to\R$ is the functional defined by
\begin{equation}\label{J}
J(f_2)={\cal L}(u(f_2),f_2)\qquad\forall\,f_2\in\R.
\end{equation}
Note that, in general the functional $J$
is not strictly convex. This implies that the solution of the optimization problem in (\ref{Jn})
is not unique and, therefore, Corollaries \ref{cor3}, \ref{cor4} cannot be applied. Nevertheless, we note that functional $J$ is stricly convex in the particular case when  $\Gamma_3=\emptyset$. Indeed, in this case problem $\cP$ consists to find an element $v\in V$ such that
\[a(u,v)=(f,v)_V\qquad\forall\, v\in V.\] Then, it is easy to see that the operator $f_2\mapsto u(f_2): L^2(\Gamma_2)\to V$ is linear and continuous and, therefore, the functional $J$ is strictly convex, which implies the unique solvability of Problem $\cQ$. We conclude that in this case Corollaries \ref{cor3} and \ref{cor4}
can be applied. We have a similar conclusion in the case when the function $g$ vanishes, i.e., $g\equiv 0$. Another important case when Corolaries \ref{cor3} and \ref{cor4} is when equality $\phi=0_V$ holds. Indeed, in this case the cost function ${\cal L}$ becomes
\begin{eqnarray*}
&&\hspace{-12mm}{\cal L}(u,f_2)=a_0\|u\|^2_{L^2(D)}+a_2\|f_2\|^2_{L^2(\Gamma_2)}
\end{eqnarray*}
for all  $u\in V$, $f_2\in L^2(\Gamma_2)$
and, it was proved in \cite{BT3} that in this case the corresponding functional $J$ defined by (\ref{J}) is strictly convex.

\medskip
{\rm (iii)} Theorem \ref{t2}  establishes a link between the solutions of  the optimal control problem $\cQ$ and the solutions of the optimal control problem of finding an element
$(u^*, f^*_2)\in {\cal V}_{ad}(\theta)$ such that}
\begin{equation}\label{3n}
{\cal L}_\theta(u^*,f^*_2)=\min_{(u,f_2)\in {\cal V}_{ad}(\theta)} {\cal L}_\theta(u,f_2).
\end{equation}
A short comparision between the optimal control problems $\cQ$ and (\ref{3n}) shows that in problem (\ref{3n}) both the state equation (and, therefore,  the set of admissible displacement pairs) and the cost functional are different to those in Problem $\cQ$. The importance of  Theorem \ref{t2} is that it provides a convergence result between the solutions of these optimal control problems which have a different structure.  In particular,  Corollary \ref{cor4} provides a continuous dependence result for the solutions of the optimal control $\cQ$ with respect to the target displacement $\phi$. This property has an important  mechanical interpretation
since it shows that, in the context of frictional antiplane shear with elastic materials, small perturbation in the target displacement field give rise to small perturbation in the corresponding optimal optimal pairs, i.e. in the optimal solution of Problem $\cQ$.

\medskip

\section*{Acknowledgement}

\indent This project has received funding from the European Union’s Horizon 2020
Research and Innovation Programme under the Marie Sklodowska-Curie
Grant Agreement No 823731 CONMECH.

\end{document}